\documentstyle[11pt]{article}

 \input amssym.def\input amssym
 
\begin{document}

\title{\bf General explicit expressions  for intertwining operators
and direct rotations of two orthogonal projections \thanks{This
research was partially supported by the National Natural Science
Foundation of China(No.11571211, 11471200), and the
Fundamental Research Funds for the Central
Universities(GK201301007). }}
\author{  \ \ Yan-Ni Dou, \ \ Wei-Juan Shi, \\
\ \ Miao-Miao Cui, \ \  Hong-Ke Du\thanks{Corresponding author:
hkdu@snnu.edu.cn. }}

\date{}
\maketitle\begin{center}
\begin{minipage}{16cm}
{\small $$College \ of \ Mathematics \ and \ Information \ Science,\
Shaanxi \ Normal \ University, $$ $$Xi'an, 710062, People's \
Republic \ of\ China. $$ }
\end{minipage}
\end{center}
\vspace{0.05cm}
\begin{center}
\begin{minipage}{16cm}

{{\bf Abstract.} In this paper, based on the block operator technique
and operator spectral theory, the general explicit expressions for
intertwining operators and direct rotations of two orthogonal
projections have been established. As a consequence, it is an improvement of Kato's result (Perturbation Theory of Linear
operators, Springer-Verlag, Berlin/Heidelberg, 1996); J. Avron, R.
Seiler and B. Simon's Theorem 2.3 (The index of a pair of
projections, J. Funct. Anal. 120(1994) 220-237) and C. Davis, W.M.
Kahan, (The rotation of eigenvectors by a perturbation, III. SIAM J.
Numer. Anal. 7(1970) 1-46). }
\endabstract

\end{minipage}\vspace{0.10cm}
\begin{minipage}{16cm}

\end{minipage}
\end{center}
\begin{center} \vspace{0.01cm}
\end{center}

\section*{1. Introduction}
\ \ \ \ Let $\mathcal{H}$ be a Hilbert space and $\mathcal{B(H)}$
the space of all bounded linear operators on $\mathcal{H}.$ An
operator $P$ is called an orthogonal projection if $P=P^*=P^2.$ Let
$\mathcal{P}$ be the set of all orthogonal projections in
$\mathcal{B(H)}.$ As well-known, orthogonal projections on a Hilbert
space are basic objects of study in operator theory (see [1-19] and
therein references). Orthogonal projections appear in various
problems and in many different areas, pure or applied. In this
paper, we will pay attention on the characterization to intertwining
operators and direct rotations of two orthogonal projections. Let the
set of all unitaries in $\mathcal{B(H)}$ be denoted by
$\mathcal{U(H)}.$ If $P$ and $Q$ are orthogonal projections and
there exists a unitary $U\in \mathcal{U(H)}$ such that
\begin{equation}UP=QU,\end{equation} then $U$ is called an outer
intertwining operator of $P$ and $Q$. The set of all outer
intertwining operators of a  pair $(P,Q)$ of orthogonal projections
is denoted by
$$\hbox{out}{\mathcal{U}}_Q(P).$$ Similarly, if
\begin{equation}PU=UQ,\end{equation} then $U$ is called an
inner intertwining operator of $P$ and $Q.$ The set of all inner
intertwining operators of a pair $(P,Q)$ of orthogonal projections
is denoted by
$$\hbox{inn}{\mathcal{U}}_Q(P).$$ Moreover, if both of
\begin{equation}PU=UQ \hbox{ and } UP=QU\end{equation} hold, then $U$
is called an intertwining operator of $P$ and $Q.$ The set of all
 intertwining operators of a pair $(P,Q)$ of orthogonal
projections is denoted by
$$\hbox{int}{\mathcal{U}}_Q(P).$$

For a pair $(P,Q)$ of orthogonal projections. A unitary $U\in
\mathcal{U(H)}$ is called a direct rotation from $P$ to $Q$ (see [1]
and [10]) if
\begin{equation}UP=QU, U^2=(Q^\perp-Q)(P^\perp-P), \hbox{Re}U\geq 0,
\end{equation} where $K^\perp=I-K$ if $K$ is an orthogonal
projection.

If $P$ and $Q$ are orthogonal projections with $\parallel
P-Q\parallel<1,$ Kato in [13] verified that there exists $U\in
\mathcal{U(H)}$ such that $PU=UQ.$  Moreover, Avron, Seiler and
Simon ([6]) proved that if $P$ and $Q$ are orthogonal projections on
$\mathcal{H}$ with $\parallel P-Q\parallel<1,$ then there exists a
unitary $U\in \mathcal{U(H)}$ with $UPU^*=Q, UQU^*=P.$
If $P$ and $Q$ are orthogonal projections have no common eigenvectors,
the mine result shown by Amrein, Sinha ([2]) implies that there exists a self-adjoint intertwining operator of $(P,Q)$.
For a pair $(P,Q)$ of orthogonal projections, we ([19]) provided a sufficient and necessary condition that there exists an
intertwining operator of $(P,Q).$ More recently, Simon ([18]) presented a more elegant proof of our previous result. In the present paper,
we will give another alternative proof of the sufficient and necessary condition for the existence of intertwining operator of $(P,Q).$
The proof is more geometrical compared with the proof in [18], and we believe the block operator technique used here has meaning in itself.

For the sake of convenience, we need some notation and terminologies. For $A\in
\mathcal{B(H)},$ the range, the null space, the spectrum, the real
part and the adjoint of $A$ denote by ${\mathcal{R}}(A),$
${\mathcal{N}}(A),$ $\sigma(A),$ $\hbox{Re}A$ and $A^*,$
respectively. $A$ is said to be positive if $(Ax,x)\geq 0$ for $x\in
\mathcal{H}.$ If $A$ is positive, then $A^\frac{1}{2}$ denotes the
positive square root of $A.$  The $A\in \mathcal{B(H)}$ is said to
be normal if $A^*A=AA^*.$ If $A$ is normal, then there exists a
spectral representation $A=\int_{\sigma(A)}\lambda dE_\lambda.$ Let
$A=U(A^*A)^\frac{1}{2}$ be the polar decomposition of $A.$  If
$\overline{{\mathcal{R}}(A)}=\mathcal{H}$ and
$\overline{{\mathcal{R}}(A^*)}=\mathcal{H},$ then $U$ in the polar
decomposition of $A$ can be chosen as a unitary. An operator $U$ is
said to be unitary if $U^*U=UU^*=I,$ where $I$ is the identity on
$\mathcal{H}.$

The following lemma is a starting point and a very useful tool in
the sequel.

{\bf Lemma 1.1.} ([11], [19]) If $W$ and $L$ are two closed
subspaces of $\mathcal{H}$ and $P$ and $Q$ denote the orthogonal
projections on $W$ and $L$, respectively, then $P$ and $Q$ have the
operator matrices
\begin{equation}P=I_1\oplus I_2\oplus 0I_3\oplus 0I_4 \oplus
I_5\oplus 0I_6\end{equation} and
\begin{equation}Q=I_1\oplus 0I_2\oplus I_3\oplus 0I_4\oplus
\left(\begin{array}{cc}
Q_0&Q_0^\frac{1}{2}(I_5-Q_0)^\frac{1}{2}D\\D^*Q_0^\frac{1}{2}
(I_5-Q_0)^\frac{1}{2}& D^*(I_5-Q_0)D\end{array}\right)\end{equation}
with respect to the space decomposition
${\mathcal{H}}=\oplus_{i=1}^6{\mathcal{H}}_i$, respectively, where
${\mathcal{H}}_1=W\cap L$, ${\mathcal{H}}_2=W\cap L^\perp$,
${\mathcal{H}}_3=W^\perp\cap L$, ${\mathcal{H}}_4=W^\perp\cap
L^\perp$, ${\mathcal{H}}_5=W\ominus({\mathcal{H}}_1\oplus
{\mathcal{H}}_2)$ and ${\mathcal{H}}_6={\mathcal{H}}\ominus(\oplus
_{j=1}^5{\mathcal{H}}_j),$ $Q_0$ is a positive contraction on
${\mathcal{H}}_5$, $0$ and $1$ are not eigenvalues of $Q_0$, and $D$
is a unitary from ${\mathcal{H}}_6$ onto ${\mathcal{H}}_5$. $I_i$ is
the identity on ${\mathcal{H}}_i$, $i=1,\ldots,6.$

{\bf Remark 1.2.} From Lemma 1.1, we will get more information
involving with geometry structure between $P$ and $Q.$ For example,

(1) Since $DD^*=I_5$ and $D^*D=I_6,$  it implies that $\dim
{\mathcal{H}}_5=\dim {\mathcal{H}}_6,$ where $\dim M$ denotes the
dimension of a subspace $M.$

(2) If $0$ ( or $1)\in \sigma(Q_0),$ then $0$ ( or $1$) is a limit
point of $\sigma(Q_0)$ and $0\notin \sigma_p(Q_0),$ where
$\sigma_p(T)$ denotes the point spectrum of $T.$ In this case, $\dim
{\mathcal{H}}_5=\dim {\mathcal{H}}_6=\infty.$ If $\dim
{\mathcal{H}}_5<\infty,$ then $0,1\notin \sigma(Q_0).$

If ${\mathcal{H}}_i=\{0\}, i=1,2,3,4,$ Halmos ([12]) called that the
pair $(P,Q)$ is in the generic position. If two orthogonal
projections are in the generic position, then
${\mathcal{H}}={\mathcal{H}}_5\oplus {\mathcal{H}}_6$ and the
operator matrices (5) and (6) of $P$ and $Q$ can be simplified as
follows
\begin{equation} P=I_5\oplus 0I_6, Q=\left(\begin{array}{cc}
Q_0&Q_0^\frac{1}{2}(I_5-Q_0)^\frac{1}{2}D\\D^*Q_0^\frac{1}{2}(I_5-Q_0)^\frac{1}{2}&
D^*(I_5-Q_0)D\end{array}\right),\end{equation} respectively. In
general, for a pair $(P,Q)$ of orthogonal projections with operator
matrices as (5) and (6), denote $\widetilde{P}$ and $\widetilde{Q}$
by  \begin{equation} \widetilde{P}=I_5\oplus 0I_6,
\widetilde{Q}=\left(\begin{array}{cc}
Q_0&Q_0^\frac{1}{2}(I_5-Q_0)^\frac{1}{2}D\\D^*Q_0^\frac{1}{2}(I_5-Q_0)^\frac{1}{2}&
D^*(I_5-Q_0)D\end{array}\right),\end{equation} the pair
$(\widetilde{P},\widetilde{Q})$ as the restriction of $(P,Q)$ on
${\mathcal{H}}_5\oplus {\mathcal{H}}_6$ is called the generic part
of $(P,Q).$

Let us give a brief outline of the contents of this paper. The
general explicit expressions for outer intertwining operators and
intertwining operators of two orthogonal projections in the generic
position are stated in Section 2. In Section 3, based on block
operator technique and spectral theory we give an alternative proof
of the sufficient and necessary condition that there exists an
intertwining operator of a pair $(P,Q).$ In view of the proof, we
get general explicit expressions of intertwining operators of a pair
$(P,Q).$ In Section 4, we provide an alternative proof of the
sufficient and necessary condition which there exists a direct
rotation of a pair $(P,Q)$ and obtain the general explicit
expressions of all direct rotations for a pair $(P,Q).$

\section*{2. General explicit expressions of intertwining operators for
a pair $(P, Q)$ in the generic position}

\qquad For outer (or inner ) intertwining operators and intertwining operators of a pair of
orthogonal projections, we have:

{\bf Theorem 2.1.} Let $P$ and $Q$ be two orthogonal projections in
the generic position and $P$ and $Q$ have operator matrix forms (7).
Then

$ \hbox{(a)}\quad \hbox{out}{\mathcal{U}}_Q(P)=
\left\{\left(\begin{array}{cc}Q_0^\frac{1}{2} &
 (I_5-Q_0)^\frac{1}{2}D
  \\D^*(I_5-Q_0)^\frac{1}{2}&-D^*Q_0^\frac{1}{2}D
   \end{array}\right)\left(\begin{array}{cc}U_0&0\\0&S_0
   \end{array}\right): U_0\in {{\mathcal{U(H}}_5)}, S_0\in
{\mathcal{U(H}}_6)\right\}.$

$ \hbox{(b)}\quad \hbox{int}{\mathcal{U}}_Q(P)=
\left\{\left(\begin{array}{cc}Q_0^\frac{1}{2} &
 (I_5-Q_0)^\frac{1}{2}D
  \\D^*(I_5-Q_0)^\frac{1}{2}&-D^*Q_0^\frac{1}{2}D
   \end{array}\right)\left(\begin{array}{cc}U_0&0\\0&D^*U_0D
   \end{array}\right):U_0\in {\mathcal{U(H}}_5),
U_0Q_0=Q_0U_0 \right\}.$

{\bf Proof.} We define an operator $W_0$ by the operator matrix
\begin{equation} W_0=\left(\begin{array}{cc}Q_0^\frac{1}{2}&(I_5-Q_0)^\frac{1}{2}D\\D^{*}(I_5-Q_0)^\frac{1}{2}&-D^{*}Q_0^\frac{1}{2}D
   \end{array}\right)\end{equation}
with the decomposition ${\mathcal{H}}=
{\mathcal{H}}_5\oplus {\mathcal{H}}_6$.
 By direct computation, $W_0$ is a unitary on ${\mathcal{H}}$ with $WP=QW$, and hence $W_0\in\hbox{int}{\mathcal{U}}_{Q}(P)\subset \hbox{out}{\mathcal{U}}_{Q}(P).$

 Note that for any $V\in\hbox{out}{\mathcal{U}}_Q(P)$, $W_0^{*}V$ is a unitary commutes with $P$,
 and for any $V\in\hbox{int}{\mathcal{U}}_Q(P)$, $W_0^{*}V$ is a unitary commutes with both $P$ and $Q$. We obtain

$$\hbox{out}{\mathcal{U}}_{Q}(P)=\{ W_0U:U\in {{\mathcal{U(H}})} \hbox{\  with\  } UP=PU  \}, $$
and
$$\hbox{int}{\mathcal{U}}_{Q}(P)=\{ W_0U:U\in {{\mathcal{U(H}})} \hbox{\  with\  } UP=PU, UQ=QU \}. $$

Let $U\in
\mathcal{U(H)}$ and $U$ has
the operator matrix
$$U=\left(\begin{array}{cc}U_{11}&U_{12}\\U_{21}&U_{22}\end{array}\right)$$
 with the decomposition ${\mathcal{H}}=
{\mathcal{H}}_5\oplus {\mathcal{H}}_6$.

If $UP=PU$, then $U_{12}=U_{21}=0$, and $U$ has the operator matrix

\begin{equation} U=\left(\begin{array}{cc}U_{11}&0\\0&U_{22}
\end{array}\right)\end{equation}
 in which $U_{11}\in {{\mathcal{U(H}}_5)}, U_{22}\in
{\mathcal{U(H}}_6)$.

Moreover, if $UP=PU$ and $UQ=QU$, then we get

\begin{equation}\left\{\begin{array}{l}U_{11}Q_{0}=Q_0U_{11},\\
U_{11}Q_{0}^\frac{1}{2}(I_5-Q_0)^\frac{1}{2}D=Q_0^\frac{1}{2}
(I_5-Q_0)^\frac{1}{2}DU_{22}.\end{array}
\right.
\end{equation}
Observing that $Q_0$,$I_5-Q_0$ on ${\mathcal{H}}_5$ are
injective from Remark 1.2 and $U_{11}$ commutes with $Q_0^\frac{1}{2},(I_5-Q_0)^\frac{1}{2}$.
It follow that $U_{11}D=DU_{22}$, and hence $U_{22}=D^{*}U_{11}D$, $U$ has the operator matrix

\begin{equation} U=\left(\begin{array}{cc}U_{11}&0\\0&D^{*}U_{11}D
\end{array}\right)\end{equation}
 in which $U_{11}\in {{\mathcal{U(H}}_5)}, U_{11}Q_{0}=Q_0U_{11}$. By (9), (10) and (12), we see that $W_0U$ in $\hbox{out}{\mathcal{U}}_Q(P)$ and $\hbox{int}{\mathcal{U}}_Q(P)$
has the operator form given in (a), (b), respectively.

The proof is completed.

{\bf Corollary 2.2.} Let a pair $(P,Q)$ of orthogonal projections be
in the generic position, and $U\in\hbox{int}{\mathcal{U}}_Q(P)$ has the operator matrix form
in Theorem 2.1. (b). Then $U$ is self-adjoint if and only if $U_0$ is self-adjoint.

{\bf Proof.} If $U$ is self-adjoint, then $U_0Q_0=Q_0U_0$ and
$Q_0^\frac{1}{2}U_0$ is self-adjoint. We get
$$U_0Q_0^\frac{1}{2}=Q_0^\frac{1}{2}U_0=U_0^*Q_0^\frac{1}{2}.$$
Hence, $(U_0-U_0^*)Q_0^\frac{1}{2}=0.$ Observing that the range of
$Q_0^\frac{1}{2}$ is dense, it follows that $U_0=U_0^*.$ This
shows that $U_0$ is self-adjoint. Conversely, it is obvious that $U$ is self-adjoint.

{\bf Remark 2.3.} (1) In Theorem 2.1. (a), the operator matrix $U\in\hbox{out}{\mathcal{U}}_Q(P)$ can be rewritten as
following, \begin{equation}U=\left(\begin{array}{cc}U_0&0\\0&D^*U_0D
   \end{array}\right)\left(\begin{array}{cc}Q_0^\frac{1}{2}
& (I_5-Q_0)^\frac{1}{2}D
  \\D^*(I_5-Q_0)^\frac{1}{2}&-D^*Q_0^\frac{1}{2}D
   \end{array}\right)\end{equation} since $U_0Q_0=Q_0U_0.$

(2)  In Theorem 2.1. (b), the operator matrix $U\in\hbox{int}{\mathcal{U}}_Q(P)$ can be rewritten as
$$U=\left(\begin{array}{cc}U_0&0\\0&S_0
   \end{array}\right)\left(\begin{array}{cc}Q_0^\frac{1}{2} &
 (I_5-Q_0)^\frac{1}{2}D
  \\D^*(I_5-Q_0)^\frac{1}{2}&-D^*Q_0^\frac{1}{2}D
   \end{array}\right),$$ where $U_0\in {\mathcal{U(H}}_5)$ and
  $S_0\in {\mathcal{U(H}}_6).$

\section*{3. General explicit expression  of intertwining operators for
two orthogonal projections}

\ \ \ \ In this section, we will devote to general explicit
expressions for intertwining operators of two orthogonal projections
if there exists an intertwining operator for the two orthogonal
projections.

Let $P$ and $Q$ be two orthogonal projections and have operator
matrices (5) and (6), respectively. For the pair $(P,Q)$ of
orthogonal projections, if the generic part of $(P,Q)$ is
$(\widetilde{P},\widetilde{Q})$ as operator matrices (8),
 then the pair
$(\widetilde{P},\widetilde{Q})$ as a pair of orthogonal projections
on ${\mathcal{H}}_5\oplus {\mathcal{H}}_6$ is in the generic
position.

The main goal in this section is to prove the following theorem.

{\bf Theorem 3.1. } Let $(P,Q)$ be a pair of orthogonal projections
with operator matrices (5) and (6), respectively. There exists a
unitary $U\in \mathcal{U(H)}$ such that $PU=UQ$ and $UP=QU$ if and
only if $\dim {\mathcal{R}}(P)\cap {\mathcal{N}}(Q)=\dim
{\mathcal{N}}(P)\cap {\mathcal{R}}(Q).$
Moreover, if $\dim {\mathcal{R}}(P)\cap {\mathcal{N}}(Q)=\dim
{\mathcal{N}}(P)\cap {\mathcal{R}}(Q),$ then
$$\begin{array}{rcl}&&\hbox{int}{\mathcal{U}}_Q(P)\\&=&
\left\{\begin{array}{l}U_1\oplus
\left(\begin{array}{cc}0&C_2\\C_3&0\end{array}\right) \oplus
U_4\oplus\left(\begin{array}{cc}Q_0^\frac{1}{2} &
 (I_5-Q_0)^\frac{1}{2}D
  \\D^*(I_5-Q_0)^\frac{1}{2}&-D^*Q_0^\frac{1}{2}D
   \end{array}\right)\left(\begin{array}{cc}U_0&0\\0&D^*U_0D
   \end{array}\right):\\U_1\in {\mathcal{U}}({\mathcal{H}}_1),
C_2\in {\mathcal{U}}({\mathcal{H}}_3,{\mathcal{H}}_2), C_3\in
{\mathcal{U}}({\mathcal{H}}_2,{\mathcal{H}}_3), U_4\in
{\mathcal{U}}({\mathcal{H}}_4), U_0\in {\mathcal{U(H}}_5),
U_0Q_0=Q_0U_0
\end{array}\right\}.\end{array}$$

{\bf Proof.} $``\Rightarrow".$ If there exists a unitary $U\in
\mathcal{U(H)}$ such that $PU=UQ$ and $UP=QU,$ then
\begin{equation}U(P-Q)=-(P-Q)U.\end{equation} Denote $A=P-Q.$ Then $A$ is a self-adjoint
contraction. So that, ${\mathcal{N}}(A),$ ${\mathcal{N}}(A-I)$ and
${\mathcal{N}}(A+I)$ are reduced subspaces of $A.$ Take
${\mathcal{H}}_0={\mathcal{H}}\ominus({\mathcal{N}}(A)\oplus
{\mathcal{N}}(A-I)\oplus{\mathcal{N}}(A+I)),$ then $A$ has the
operator matrix
\begin{equation}A=\left(\begin{array}{cccc}0&0&0&0\\0&I_1&0&0\\0&0&
-I_{-1}&0\\0&0&0&A_0\end{array}\right)\end{equation} with respect to
the decomposition ${\mathcal{H}}={\mathcal{N}}(A)\oplus
{\mathcal{N}}(A-I)\oplus{\mathcal{N}}(A+I)\oplus {\mathcal{H}}_0,$
where $I_1$ is the identity on ${\mathcal{N}}(A-I),$ $I_{-1}$ is the
identity on ${\mathcal{N}}(A-I),$ $I_0$ is the identity on
${\mathcal{H}}_0.$

It is clear that $A_0,$ $A_0-I_0$ and $A_0+I_0$ as operators on
${\mathcal{H}}_0$ are injective and dense.

 If ${U}$ has the operator matrix
 \begin{equation}U=\left(\begin{array}{cccc}U_{11}&U_{12}&U_{13}&U_{14}\\
 U_{21}&U_{22}&U_{23}&U_{24}\\U_{31}&U_{32}&U_{33}&U_{34}
 \\U_{41}&U_{42}&U_{43}&U_{44}\end{array}\right)\end{equation}
  with respect to the
decomposition ${\mathcal{H}}={\mathcal{N}}(A)\oplus
{\mathcal{N}}(A-I)\oplus{\mathcal{N}}(A+I)\oplus {\mathcal{H}}_0,$
then from (14), we get $UA=-AU.$ Moreover, by (15) and (16), we
obtain
\begin{equation}\left(\begin{array}{cccc}
0&U_{12}&-U_{13}&U_{14}A_0\\
 0&U_{22}&-U_{23}&U_{24}A_0\\0&U_{32}&-U_{33}&U_{34}A_0
 \\0&U_{42}&-U_{43}&U_{44}A_0\end{array}\right)=
 \left(\begin{array}{cccc}0&0&0&0\\
 -U_{21}&-U_{22}&-U_{23}&-U_{24}\\U_{31}&U_{32}&U_{33}&U_{34}
 \\-A_0U_{41}&-A_0U_{42}&-A_0U_{43}&-A_0U_{44}\end{array}\right).
 \end{equation}
 Comparing two sides of (17) and observing that
$A_0,$ $A_0-I_0$ and $A_0+I_0$ are injective and dense, it is
derived that $U_{12}=0,  U_{13}=0, U_{14}=0, U_{21}=0, U_{22}=0,
U_{24}=0, U_{31}=0, U_{33}=0, U_{34}=0, U_{41}=0, U_{42}=0,
U_{43}=0.$ Therefore,
\begin{equation}U=\left(\begin{array}{cccc}U_{11}&0&0&0\\
 0&0&U_{23}&0\\0&U_{32}&0&0
 \\0&0&0&U_{44}\end{array}\right)\end{equation} This shows that
 $U_{11}$ is a unitary on ${\mathcal{N}}(A)={\mathcal{H}}_1\oplus
 {\mathcal{H}}_4,$
 $\left(\begin{array}{cc}0&U_{23}\\U_{32}&0\end{array}\right)$ is a
 unitary on ${\mathcal{N}}(A-I)\oplus
 {\mathcal{N}}(A+I)=({\mathcal{R}}(P)\cap {\mathcal{N}}(Q))\oplus
 ({\mathcal{N}}(P)\cap {\mathcal{R}}(Q))$ and $U_{44}$ is a unitary on
 ${\mathcal{H}}_0={\mathcal{H}}_5\oplus {\mathcal{H}}_6.$

Observing that $U_{11}P_{{\mathcal{H}}_1\oplus
 {\mathcal{H}}_4}=Q_{{\mathcal{H}}_1\oplus
 {\mathcal{H}}_4}U_{11}$ and $U_{11},$ $P_{{\mathcal{H}}_1\oplus
 {\mathcal{H}}_4}$ and $Q_{{\mathcal{H}}_1\oplus
 {\mathcal{H}}_4}$ have operator matrices
 $$U_{11}=\left(\begin{array}{cc}U_{11}^{11}&U_{11}^{12}\\
 U_{11}^{21}&U_{11}^{22}\end{array}\right), P_{{\mathcal{H}}_1\oplus
 {\mathcal{H}}_4}=Q_{{\mathcal{H}}_1\oplus
 {\mathcal{H}}_4}=\left(\begin{array}{cc}
 I_1&0\\0&0\end{array}\right)$$ with respect to the decomposition
 ${\mathcal{H}}_1\oplus
 {\mathcal{H}}_4,$ respectively, from $UP=QU$ we get
 $$\left(\begin{array}{cc}U_{11}^{11}&0\\U_{11}^{21}&0
 \end{array}\right)=\left(\begin{array}{cc}U_{11}^{11}&U_{11}^{12}
 \\0&0\end{array}\right).$$ Hence, $U_{11}^{12}=0$ and
 $U_{11}^{21}=0.$ Therefore, $U_{11}^{11}$ and $U_{11}^{22}$ are
 unitaries on ${\mathcal{H}}_1$ and ${\mathcal{H}}_4,$ respectively.

Observing that $U_{{\mathcal{H}}_2\oplus
 {\mathcal{H}}_3}=\left(\begin{array}{cc}0&U_{23}\\U_{32}&0
\end{array}\right)$ and $U_{{\mathcal{H}}_2\oplus
 {\mathcal{H}}_3}$ is a unitary, we have
 $$\left(\begin{array}{cc}0&U_{23}\\U_{32}&0
\end{array}\right)\left(\begin{array}{cc}0&U_{32}^*\\U_{23}^*&0
\end{array}\right)=\left(\begin{array}{cc}U_{23}U_{23}^*&0\\
0&U_{32}U_{32}^*
\end{array}\right)=\left(\begin{array}{cc}I_2&0\\0&I_3
\end{array}\right)$$ and $$\left(\begin{array}{cc}0&U_{32}^*\\U_{23}^*&0
\end{array}\right)\left(\begin{array}{cc}0&U_{23}\\U_{32}&0
\end{array}\right)=\left(\begin{array}{cc}U_{32}^*U_{32}&0\\
0&U_{23}^*U_{23}
\end{array}\right)=\left(\begin{array}{cc}I_2&0\\0&I_3
\end{array}\right).$$ Thus $U_{32}^*U_{32}=I_2$ and
$U_{23}^*U_{23}=I_3.$ It implies that $\dim{\mathcal{H}}_2=\dim
 {\mathcal{H}}_3$ and $U_{23}$ is a unitary from ${\mathcal{H}}_3$
 onto ${\mathcal{H}}_2.$ Similarly, $U_{32}$ is a unitary from
  ${\mathcal{H}}_2$
 onto ${\mathcal{H}}_3.$

 Next, from $UP=QU$ and $UQ=PU,$ we have
 $$U_{44}P_{{\mathcal{H}}_5\oplus{\mathcal{H}}_6}=
 Q_{{\mathcal{H}}_5\oplus{\mathcal{H}}_6}
U_{44} \hbox{ and }U_{44}Q_{{\mathcal{H}}_5\oplus{\mathcal{H}}_6}=
P_{{\mathcal{H}}_5\oplus{\mathcal{H}}_6} U_{44}.$$ By Theorem 2.1,
$$U_{44}=\left(\begin{array}{cc}Q_0^\frac{1}{2} &
 (I_5-Q_0)^\frac{1}{2}D
  \\D^*(I_5-Q_0)^\frac{1}{2}&-D^*Q_0^\frac{1}{2}D
   \end{array}\right)\left(\begin{array}{cc}U_0&0\\0&D^*U_0D
   \end{array}\right),$$ where $U_0$ with $U_0Q_0=Q_0U_0 $ is a unitary
on ${\mathcal{H}}_5.$

$``\Leftarrow". $ If $\dim {\mathcal{R}}(P)\cap
{\mathcal{N}}(Q)=\dim {\mathcal{N}}(P)\cap {\mathcal{R}}(Q),$ we can
choose a unitary $C_2$ from ${\mathcal{H}}_3$ onto ${\mathcal{H}}_2$
and a unitary $C_3$ from ${\mathcal{H}}_2$ onto ${\mathcal{H}}_3.$
Define an operator $$U=U_1\oplus
\left(\begin{array}{cc}0&C_2\\C_3&0\end{array}\right)\oplus
U_4\oplus \left(\begin{array}{cc}Q_0^\frac{1}{2}U_0 &
 (I_5-Q_0)^\frac{1}{2}U_0D
  \\D^*(I_5-Q_0)^\frac{1}{2}U_0&-D^*Q_0^\frac{1}{2}U_0D
   \end{array}\right),$$ where $U_1$ is a unitary on
${\mathcal{H}}_1,$ $C_2$ is a unitary from ${\mathcal{H}}_3$ onto
${\mathcal{H}}_2,$ $C_3$ is a unitary from ${\mathcal{H}}_2$ onto
${\mathcal{H}}_3,$ $U_4$ is a unitary on ${\mathcal{H}}_4$ and $U_0$
is a unitary on ${\mathcal{H}}_5$ with $Q_0U_0=U_0Q_0,$ by
directly checking, $U$ is a unitary on $\mathcal{H}$ and $UP=QU$ and
$ UQ=PU.$

From the proof above, we have
$$\begin{array}{rcl}&&\hbox{int}{\mathcal{U}}_Q(P)\\&=&
\left\{\begin{array}{l}U_1\oplus
\left(\begin{array}{cc}0&C_2\\C_3&0\end{array}\right) \oplus
U_4\oplus\left(\begin{array}{cc}Q_0^\frac{1}{2} &
 (I_5-Q_0)^\frac{1}{2}D
  \\D^*(I_5-Q_0)^\frac{1}{2}&-D^*Q_0^\frac{1}{2}D
   \end{array}\right)\left(\begin{array}{cc}U_0&0\\0&D^*U_0D
   \end{array}\right):\\U_1\in {\mathcal{U}}({\mathcal{H}}_1),
C_2\in {\mathcal{U}}({\mathcal{H}}_3,{\mathcal{H}}_2), C_3\in
{\mathcal{U}}({\mathcal{H}}_2,{\mathcal{H}}_3), U_4\in
{\mathcal{U}}({\mathcal{H}}_4), U_0\in {\mathcal{U(H}}_5),
U_0Q_0=Q_0U_0
\end{array}\right\}.\end{array}$$

{\bf Remark 3.2.}  Let $(P,Q)$ be a pair of orthogonal projections.
From the proof of Theorem 3.1, if $\dim {\mathcal{R}}(P)\cap
{\mathcal{N}}(Q)=\dim {\mathcal{N}}(P)\cap {\mathcal{R}}(Q),$ then
the intertwining operator of the pair of orthogonal projections is
not unique. Moreover, it can be chose as a self-adjoint unitary.
Even though the intertwining operator can be chose as a self-adjoint
operator, it is also not unique by Corollary 2.2.

As a consequence, we give an alternative proof of Theorem 2.2 in
[15].

{\bf Corollary 3.3.} (Theorem 2.2 in [15]) Let $L$ and $M$ be
subspaces of $\mathcal{H}.$ If $P$ and $Q$ are orthogonal
projections on $L$ and $M$, respectively, then there exists a
unitary operator $U\in \mathcal{B(H)}$ such that
\begin{equation}PQP=UQPQU^*.\end{equation}

{\bf Proof.} Let $P$ and $Q$ have operator matrices (5) and (6),
respectively. Then $$PQP=I_1\oplus 0I_2\oplus 0I_3\oplus 0I_4\oplus
Q_0\oplus 0I_6$$ and
$$QPQ=I_1\oplus 0I_2\oplus 0I_3\oplus
0I_4\oplus\left(\begin{array}{cc}Q_0^2&Q_0^\frac{3}{2}
(I_5-Q_0)^\frac{1}{2}D\\D^*Q_0^\frac{3}{2}
(I_5-Q_0)^\frac{1}{2}&D^*Q_0(I_5-Q_0)D\end{array}\right).$$ Denote
the generic part $(\widetilde{P},\widetilde{Q})$ of $(P,Q)$ as the
operator matrices (8). We get
$$\widetilde{P}\widetilde{Q}\widetilde{P}=\left(\begin{array}{cc}
Q_0&0\\0&0\end{array}\right),\widetilde{Q}\widetilde{P}\widetilde{Q}=\left(\begin{array}{cc}Q_0^2&Q_0^\frac{3}{2}
(I_5-Q_0)^\frac{1}{2}D\\D^*Q_0^\frac{3}{2}
(I_5-Q_0)^\frac{1}{2}&D^*Q_0(I_5-Q_0)D\end{array}\right).$$ By
Theorem 2.1, there exists a unitary $\widetilde{U}$ on
${\mathcal{H}}_5\oplus {\mathcal{H}}_6$ such that
\begin{equation}\widetilde{U}\widetilde{P}\widetilde{U}^*=
\widetilde{Q},
\widetilde{U}\widetilde{Q}\widetilde{U}^*=\widetilde{P}.
\end{equation}
In this case,
$$\widetilde{P}\widetilde{Q}\widetilde{P}=
\widetilde{U}\widetilde{Q}\widetilde{U}^*\widetilde{U}\widetilde{P}
\widetilde{U}^*\widetilde{U}\widetilde{Q}
\widetilde{U}^*=\widetilde{U}\widetilde{Q}\widetilde{P}
\widetilde{Q} \widetilde{U}^*.$$ Furthermore, define $U$ by
\begin{equation}U=\oplus_{i=1}^4I_i\oplus\widetilde{U},\end{equation}
 where $I_i$ are identities on ${\mathcal{H}}_i,$ $1\leq i\leq 4.$

Evidently, $U$ is a unitary, and $Q=UPU^*$ and $P=UQU^*.$ Hence,
$$\begin{array}{rcl}PQP&=&I_1\oplus 0I_2\oplus 0I_3\oplus 0I_4\oplus
Q_0\oplus 0I_6\\&=&I_1\oplus 0I_2\oplus 0I_3\oplus 0I_4\oplus
\widetilde{P}\widetilde{Q}\widetilde{P}\\&=&I_1\oplus 0I_2\oplus
0I_3\oplus 0I_4\oplus \widetilde{U}\widetilde{Q}\widetilde{P}
\widetilde{Q} \widetilde{U}^* \\&=&UQPQU^*.\end{array}$$

{\bf Remark 3.4.} (1) By Theorem 2.1 and Theorem 3.1, a unitary
satisfying (19) is not unique.

(2) $U$ in Corollary 3.3 can by chose as a self-adjoint unitary.
Even so this choice is not unique by Corollary 2.2.

\section*{4. General explicit expression of direct rotations on
 a pair of orthogonal projections}

\qquad The concept of a direct rotation of a pair on orthogonal projections
due to Davis (see [10]).

{\bf Definition 4.1.} (Definition 2.9 in [1], Definition 3.1 in
[10]) Let $(P,Q)$ be a pair of orthogonal projections. A unitary
$S\in \mathcal{U(H)}$ is called a direct rotation from $P$ to $Q$
(see [10]) if
$$SP=QS, S^2=(Q^\perp-Q)(P^\perp-P), \hbox{Re}S\geq 0.$$

For a pair $(P,Q)$ of orthogonal projections, denote the set of all
direct rotations from $P$ to $Q$ by
$${\mathcal{S}}_Q(P)=\{S\in {\mathcal{U(H)}}:SP=QS,
S^2=(Q^\perp-Q)(P^\perp-P), \hbox{Re}S\geq 0\}.$$

{\bf Lemma 4.2.} (Proposition 3.1 in [10]) If a pair $(P,Q)$ of
orthogonal projections is in the generic position, then there exists
a unique unitary operator $S$ such that \begin{equation}SP=QS,
S^2=(Q^\perp-Q)(P^\perp-P), \hbox{Re}S\geq 0.\end{equation}
Moreover, if $P$ and $Q$ have the operator matrices (7), then
\begin{equation}S=\left(\begin{array}{cc}S_{11}&S_{12}\\S_{21}&S_{22}
\end{array}\right)=\left(\begin{array}{cc}Q_0^\frac{1}{2} &
 -(I_5-Q_0)^\frac{1}{2}D
  \\D^*(I_5-Q_0)^\frac{1}{2}&D^*Q_0^\frac{1}{2}D
   \end{array}\right).\end{equation}

{\bf Proof.} If there exists a unitary operator $S$ satisfying (22),
then from $SP=QS$ we get
\begin{equation}SP^\perp=Q^\perp S, S^*Q=PS^*, S^*Q^\perp=P^\perp
S^*. \end{equation} Hence, from $S^2=(Q^\perp-Q)(P^\perp-P),$ we
obtain
\begin{equation}S=S^*(Q^\perp-Q)(P^\perp-P)=
(P^\perp-P)S^*(P^\perp-P).\end{equation} Let $P,$ $Q$ and $S$ have
operator matrices
$$P=\left(\begin{array}{cc}I&0\\0&0\end{array}\right),
Q=\left(\begin{array}{cc}
Q_0&Q_0^\frac{1}{2}(I_5-Q_0)^\frac{1}{2}D\\D^*Q_0^\frac{1}{2}
(I_5-Q_0)^\frac{1}{2}& D^*(I_5-Q_0)D\end{array}\right),
S=\left(\begin{array}{cc}S_{11}&S_{12}\\S_{21}&S_{22}
\end{array}\right)$$ with respect to the decomposition
${\mathcal{H}}={\mathcal{R}}(P)\oplus{\mathcal{N}}(P),$
respectively. From (25),
$$ \left(\begin{array}{cc}S_{11}&S_{12}\\S_{21}&S_{22}
\end{array}\right)=\left(\begin{array}{cc}S_{11}^*&-S_{21}^*\\
-S_{12}^*&S_{22}^*
\end{array}\right).$$ So that,
\begin{equation}\left\{\begin{array}{l}S_{11}=S_{11}^*,\\S_{12}=-S_{21}^*,
\\S_{21}=-S_{12}^*,\\S_{22}=S_{22}^*.\end{array}\right.\end{equation}
Hence, $$S=\left(\begin{array}{cc}S_{11}&S_{12}\\-S_{12}^*&S_{22}
\end{array}\right)=\left(\begin{array}{cc}S_{11}&-S_{21}^*\\S_{21}&S_{22}
\end{array}\right).$$

 Moreover,
$$\hbox{Re}S=\frac{1}{2}(S+S^*)=
\left(\begin{array}{cc}S_{11}&0\\0&S_{22}
\end{array}\right)\geq 0.$$
In general, by Theorem 2.1, there exist two unitaries $U_0,V_0\in
{\mathcal{U}}({\mathcal{R}}(P))$ such that
$$\left(\begin{array}{cc}S_{11}&S_{12}\\S_{21}&S_{22}
\end{array}\right)=\left(\begin{array}{cc}Q_0^\frac{1}{2} &
 (I_5-Q_0)^\frac{1}{2}D
  \\D^*(I_5-Q_0)^\frac{1}{2}&-D^*Q_0^\frac{1}{2}D
   \end{array}\right)\left(\begin{array}{cc}U_0&0\\0&D^*V_0D
   \end{array}\right).$$ So that, \begin{equation}
   Q_0^\frac{1}{2}U_0=S_{11}\geq
   0\end{equation} and \begin{equation}
   -D^*Q_0^\frac{1}{2}V_0D=S_{22}\geq 0.\end{equation} Since $Q_0$
is injective, by (27) and (28), it is clear that
$U_0=I_{{\mathcal{R}}(P)},$  $V_0=-I_{{\mathcal{R}}(P)}.$ Therefore,
$$S=\left(\begin{array}{cc}S_{11}&S_{12}\\S_{21}&S_{22}
\end{array}\right)=\left(\begin{array}{cc}Q_0^\frac{1}{2} &
 -(I_5-Q_0)^\frac{1}{2}D
  \\D^*(I_5-Q_0)^\frac{1}{2}&D^*Q_0^\frac{1}{2}D
   \end{array}\right),$$ it is uniquely determined.

   If $$S=\left(\begin{array}{cc}Q_0^\frac{1}{2} &
 -(I_5-Q_0)^\frac{1}{2}D
  \\D^*(I_5-Q_0)^\frac{1}{2}&D^*Q_0^\frac{1}{2}D
   \end{array}\right),$$ by directly checking, $S$ satisfies (22).
It is the direct rotation from $P$ to $Q.$

{\bf Theorem 4.3.} (Proposition 3.2 in [10]) For a pair $(P,Q)$ of
orthogonal projections, there exists a direct rotation $S$ from $P$
to $Q$ which satisfies (22) if and only if $\dim
{\mathcal{R}}(P)\cap {\mathcal{N}}(Q)=\dim {\mathcal{N}}(P)\cap
{\mathcal{R}}(Q).$
Moreover, if $P$ and $Q$ with $\dim {\mathcal{R}}(P)\cap
{\mathcal{N}}(Q)=\dim {\mathcal{N}}(P)\cap {\mathcal{R}}(Q)$ have
the operator matrices (5) and (6) with respect to the space
decomposition ${\mathcal{H}}=\oplus_{i=1}^6{\mathcal{H}}_i,$ then
\begin{equation}S=I_1\oplus
\left(\begin{array}{cc}0&C\\-C^*&0\end{array}\right)\oplus I_4\oplus
\left(\begin{array}{cc}Q_0^\frac{1}{2} &
 -(I_5-Q_0)^\frac{1}{2}D
  \\D^*(I_5-Q_0)^\frac{1}{2}&D^*Q_0^\frac{1}{2}D
   \end{array}\right),\end{equation} where $C$ is an arbitrary unitary from
 ${\mathcal{H}}_3$ onto ${\mathcal{H}}_2.$

 {\bf Proof.} Denote
${\mathcal{K}}_1=({\mathcal{R}}(P)\cap{\mathcal{R}}(Q))
\oplus({\mathcal{N}}(P)\cap{\mathcal{N}}(Q)),$
${\mathcal{K}}_2=({\mathcal{R}}(P)\cap{\mathcal{N}}(Q))
\oplus({\mathcal{N}}(P)\cap{\mathcal{R}}(Q))$ and $
{\mathcal{K}}_3={\mathcal{H}}\ominus({\mathcal{K}}_1\oplus
{\mathcal{K}}_2) .$

For $x_1\in {\mathcal{K}}_1,$  we get
\begin{equation}S^2x_1=(Q^\perp-Q)(P^\perp-P)x_1=x_1.\end{equation}
  From (30),
we obtain
 that $(S^2-I)x_1=(S+I)(S-I)x_1=0.$
Moreover, observing that $S+I$ is invertible since $\hbox{Re}S\geq
0,$ we get $(S-I)x_1=0.$ Hence,
 $$Sx_1=x_1.$$ This shows that ${\mathcal{K}}_1$ is
 a reduced subspace under $S$ and
 $S\mid_{{\mathcal{K}}_1}$ is the identity on ${\mathcal{K}}_1.$

 For any $y\in {\mathcal{K}}_2,$ denote $y=y_1+y_2,$ where
 $y_1\in{\mathcal{R}}(P)\cap
{\mathcal{N}}(Q)$ and $ y_2\in
{\mathcal{N}}(P)\cap{\mathcal{R}}(Q),$ we shall show that $Sy\in
{\mathcal{K}}_2.$ Observing that $Sy_1=SPy_1=QSy_1\in
{\mathcal{R}}(Q)$ and $Sy_1=SQ^\perp y_1=P^\perp Sy_1\in
{\mathcal{N}}(P),$ we have $$Sy_1\in {\mathcal{N}}(P)\cap
{\mathcal{R}}(Q).$$ Similarly,  $$Sy_2\in {\mathcal{R}}(P)\cap
{\mathcal{N}}(Q).$$ Hence, $$Sy=Sy_1+Sy_2\in {\mathcal{K}}_2.$$ This
shows that ${\mathcal{K}}_2$ is an invariant subspace of $S.$ In
this case, $S$ has the operator matrix
\begin{equation}S=\left(\begin{array}{ccc}I_{{\mathcal{K}}_1}&0&0\\
0&S_{22}&S_{23}\\
0&0&S_{33}\end{array}\right)\end{equation} with respect to the
decomposition ${\mathcal{H}}=\oplus_{i=1}^3{\mathcal{K}}_i.$
Furthermore, since
$$S^2=\left(\begin{array}{ccc}I_{{\mathcal{K}}_1}&0&0\\0&S_{22}^2&
S_{22}S_{23}+S_{23}S_{33}\\
0&0&S_{33}^2\end{array}\right),$$ if $y=y_1+y_2,$ where $y_1\in
{\mathcal{R}}(P)\cap{\mathcal{N}}(Q)$ and $y_2\in
{\mathcal{N}}(P)\cap{\mathcal{R}}(Q),$ we get
$S^2y=S^2(y_1+y_2)=-y_1-y_2=-y.$ So that $S_{22}^2y=-y.$ This means
that
\begin{equation}S_{22}^2=-I_{{\mathcal{K}}_2}.\end{equation}
It implies that $S_{22}$ is an invertible operator on
${\mathcal{K}}_2.$ Furthermore,$$S^*S=
\left(\begin{array}{ccc}I_{{\mathcal{K}}_1}&0&0\\0&S_{22}^*S_{22}
&S_{22}^*S_{23}\\
0&S_{23}^*S_{22}&S_{23}^*S_{23}+S_{33}^*S_{33}\end{array}\right)=
\left(\begin{array}{ccc}I_{{\mathcal{K}}_1}&0&0\\0&
I_{{\mathcal{K}}_2}&0
\\0&0&I_{{\mathcal{K}}_3}\end{array}\right).$$
It follows that $S_{23}^*S_{22}=0.$ From (32),
${\mathcal{R}}(S_{22})={\mathcal{K}}_2,$ it is derived that
$S_{23}=0.$

So that, the operator matrix form (31) can be changed as follows
\begin{equation}S=
\left(\begin{array}{ccc}I_{{\mathcal{K}}_1}&0&0\\
0&S_{22}&0\\
0&0&S_{33}\end{array}\right).\end{equation} Here, $S_{22}$ and
$S_{33}$ are unitaries on ${\mathcal{K}}_2$ and ${\mathcal{K}}_3$,
respectively.

If $P$ and $Q$ have the operator matrices (5) and (6), then it is
obvious that $S_{22}$ as a unitary on
${\mathcal{K}}_2={\mathcal{H}}_2\oplus {\mathcal{H}}_3$ has the
operator matrix form
$$S_{22}=\left(\begin{array}{cc}0&C\\-C^*&0\end{array}\right)$$
with respect to the decomposition
${\mathcal{K}}_2={\mathcal{H}}_2\oplus {\mathcal{H}}_3,$ where $C$
is an arbitrary unitary from ${\mathcal{H}}_3$ onto
${\mathcal{H}}_2.$ Explicitly, there exists a unitary such as $C$
above if and only if $\dim {\mathcal{R}}(P)\cap{\mathcal{N}}(Q)=
\dim {\mathcal{N}}(P)\cap{\mathcal{R}}(Q).$ By Lemma 4.2, $S_{33}$
has the operator matrix form
$$S_{33}=\left(\begin{array}{cc}Q_0^\frac{1}{2} &
 -(I_5-Q_0)^\frac{1}{2}D
  \\D^*(I_5-Q_0)^\frac{1}{2}&D^*Q_0^\frac{1}{2}D
   \end{array}\right).$$ It is uniquely determined.
So that, \begin{equation}S=I_1\oplus
\left(\begin{array}{cc}0&C\\-C^*&0\end{array}\right)\oplus I_4\oplus
\left(\begin{array}{cc}Q_0^\frac{1}{2} &
 -(I_5-Q_0)^\frac{1}{2}D
  \\D^*(I_5-Q_0)^\frac{1}{2}&D^*Q_0^\frac{1}{2}D
   \end{array}\right).\end{equation}

Conversely, if $\dim {\mathcal{R}}(P)\cap{\mathcal{N}}(Q)= \dim
{\mathcal{N}}(P)\cap{\mathcal{R}}(Q),$ for any unitary $C$ from
${\mathcal{H}}_3$ onto ${\mathcal{H}}_2,$ define an operator $S$ by
the form (34), then to directly test the operator $S$ is a unitary
which satisfies (22). That is, $S$ is a direct rotation of the pair
$(P,Q)$ from $P$ to $Q.$

{\bf Remark 4.4.} (1) There exists a unique unitary $S$ satisfying
(22) if and only if $\dim {\mathcal{R}}(P)\cap{\mathcal{N}}(Q)= \dim
{\mathcal{N}}(P)\cap{\mathcal{R}}(Q)=0.$

(2) For a pair $(P,Q)$ of orthogonal projections, if $\dim
{\mathcal{R}}(P)\cap{\mathcal{N}}(Q)= \dim
{\mathcal{N}}(P)\cap{\mathcal{R}}(Q)\neq 0,$ then the direct
rotation from $P$ to $Q$ is not unique. The general expression of
direct rotations $S$ from $P$ to $Q$ has the form (29), where $C$ can
be chose over all unitaries from ${\mathcal{H}}_3$ onto
${\mathcal{H}}_2.$

(3) For a pair $(P,Q)$ of orthogonal projections with $\dim
{\mathcal{R}}(P)\cap{\mathcal{N}}(Q)= \dim
{\mathcal{N}}(P)\cap{\mathcal{R}}(Q),$ if the set of all direct
rotations from $P$ to $Q$ is denoted by ${\mathcal{S}}_Q(P),$ then
$${\mathcal{S}}_Q(P)=\left\{I_1\oplus
\left(\begin{array}{cc}0&C\\-C^*&0\end{array}\right)\oplus I_4\oplus
\left(\begin{array}{cc}Q_0^\frac{1}{2} &
 -(I_5-Q_0)^\frac{1}{2}D
  \\D^*(I_5-Q_0)^\frac{1}{2}&D^*Q_0^\frac{1}{2}D
   \end{array}\right): C \in {\mathcal{U}}({\mathcal{H}}_3,
{\mathcal{H}}_2)\right\}.$$

(4) It is interesting that if a pair $(P,Q)$ of orthogonal
projections with $\dim {\mathcal{R}}(P)\cap{\mathcal{N}}(Q)= \dim
{\mathcal{N}}(P)\cap{\mathcal{R}}(Q)$ and ${\mathcal{H}}_5\neq
\{0\},$ then $$\hbox{int}_Q(P)\cap {\mathcal{S}}_Q(P)=\emptyset.$$

As the end, we will give an alternative proof of the extremal
property in regard to the direct rotation which is due to Davis (
see [1],[10]). The proof used block operator matrices and spectral
theory may give us some inspiration in the further study.

{\bf Theorem 4.5.} Let the pair $(P,Q)$ of orthogonal projections be
in the generic position. The direct rotation $U$ from $P$ to $Q$ has
the extremal property
$$\parallel U-I\parallel=\inf\{\parallel\widetilde{U}-I\parallel:
\widetilde{U}\in {\mathcal{U(H)}},
P=\widetilde{U}^*Q\widetilde{U}\}.$$

{\bf Proof.} Assume that $P$ and $Q$ are in the generic position and
have the operator matrix (7). From Lemma 4.2 and (23), the direct
rotation $U$ from $P$ to $Q$ is unique and
$$U=\left(\begin{array}{cc}Q_0^\frac{1}{2} &
 -(I_5-Q_0)^\frac{1}{2}D
  \\D^*(I_5-Q_0)^\frac{1}{2}&D^*Q_0^\frac{1}{2}D
   \end{array}\right).$$ Hence,
$$\begin{array}{rcl}&&\parallel U-I\parallel^2\\&=&\parallel
\left(\begin{array}{cc}Q_0^\frac{1}{2}-I_5 &
 -(I_5-Q_0)^\frac{1}{2}D
  \\D^*(I_5-Q_0)^\frac{1}{2}&D^*(Q_0^\frac{1}{2}-I_5)D
   \end{array}\right)\parallel^2\\&=&\parallel
   \left(\begin{array}{cc}Q_0^\frac{1}{2}-I_5 &
 -(I_5-Q_0)^\frac{1}{2}D
  \\D^*(I_5-Q_0)^\frac{1}{2}&D^*(Q_0^\frac{1}{2}-I_5)D
   \end{array}\right)\left(\begin{array}{cc}Q_0^\frac{1}{2}-I_5 &
 (I_5-Q_0)^\frac{1}{2}D
  \\-D^*(I_5-Q_0)^\frac{1}{2}&D^*(Q_0^\frac{1}{2}-I_5)D
   \end{array}\right)\parallel\\&=&\parallel
\left(\begin{array}{cc}2(I_5-Q_0^\frac{1}{2}) &
 0\\0&2D^*(I_5-Q_0^\frac{1}{2})D
   \end{array}\right)\parallel\\&=&2\parallel I_5-Q_0^\frac{1}{2}
   \parallel. \end{array}$$

   If $\lambda_0=\min\{\lambda:\lambda\in \sigma(Q_0)\},$ then $$\parallel I_5-Q_0^\frac{1}{2}
   \parallel=1-\lambda_0^\frac{1}{2}.$$
Thus  $$\parallel
U-I\parallel=\sqrt{2(1-\lambda_0^\frac{1}{2})}.$$

By Theorem 2.1, if $Q= \widetilde{U}P\widetilde{U}^*,$  then we have
$$\widetilde{U}=\left(\begin{array}{cc}Q_0^\frac{1}{2} &
 (I_5-Q_0)^\frac{1}{2}D
  \\D^*(I_5-Q_0)^\frac{1}{2}&-D^*Q_0^\frac{1}{2}D
   \end{array}\right)\left(\begin{array}{cc}V_0&0\\0&DS_0D^*
   \end{array}\right),$$ where
 $V_0,S_0\in {\mathcal{U(H}}_5).$ In this case,
$$\begin{array}{rcl}&&\parallel \widetilde{U}-I\parallel^2\\&=&
\parallel
\left(\begin{array}{cc}Q_0^\frac{1}{2}-V_0^* &
 (I_5-Q_0)^\frac{1}{2}D
  \\D^*(I_5-Q_0)^\frac{1}{2}&-D^*(Q_0^\frac{1}{2}+S_0^*)D
   \end{array}\right)\parallel^2\\&=&\parallel
   \left(\begin{array}{cc}Q_0^\frac{1}{2}-V_0 &
 (I_5-Q_0)^\frac{1}{2}D
  \\D^*(I_5-Q_0)^\frac{1}{2}&-D^*(Q_0^\frac{1}{2}+S_0)D
   \end{array}\right)\left(\begin{array}{cc}Q_0^\frac{1}{2}-V_0^* &
 (I_5-Q_0)^\frac{1}{2}D
  \\D^*(I_5-Q_0)^\frac{1}{2}&-D^*(Q_0^\frac{1}{2}+S_0^*)D
   \end{array}\right)\parallel\\&=&\parallel
\left(\begin{array}{cc}2I_5-( V_0Q_0^\frac{1}{2}+
Q_0^\frac{1}{2}V_0^*) & *\\*
 * & D^*(2I_5+Q_0^\frac{1}{2}S_0^*+S_0Q_0^\frac{1}{2})D
   \end{array}\right)\parallel\\&\geq &\max\{\parallel
2I_5-(  Q_0^\frac{1}{2}V_0^*+V_0Q_0^\frac{1}{2}) \parallel,
\parallel 2I_5+Q_0^\frac{1}{2}S_0^*+S_0Q_0^\frac{1}{2}
   \parallel\}. \end{array}$$

Without loss of generality, we can assume that $\lambda_0\in
\sigma_p(Q_0).$ Take a unit vector $x_{\lambda_0}$ such that
$Q_0x_{\lambda_0}=\lambda_0x_{\lambda_0}.$ We get
$$\begin{array}{rcl}\parallel 2I_5-(Q_0^\frac{1}{2}V_0^*+
V_0Q_0^\frac{1}{2})\parallel&\geq &((2I_5-(V_0Q_0^\frac{1}{2}+
Q_0^\frac{1}{2}V_0^*))x_{\lambda_0},x_{\lambda_0})\\&=&
2-\lambda_0^\frac{1}{2}((V_0^*+V_0)x_{\lambda_0},x_{\lambda_0})\\&
\geq & 2(1-\lambda_0^\frac{1}{2}).
\end{array}$$
Similarly, $$\parallel 2I_5+Q_0^\frac{1}{2}S_0^*+S_0Q_0^\frac{1}{2}
   \parallel\geq 2(1-\lambda_0^\frac{1}{2}).$$ So that, $\parallel \widetilde{U}-I\parallel\geq
\sqrt{2(1-\lambda_0^\frac{1}{2})}.$ Hence, $\parallel
\widetilde{U}-I\parallel\geq \parallel U-I\parallel.$ Thus
$$\parallel U-I\parallel=\inf\{\parallel\widetilde{U}-I\parallel:
\widetilde{U}\in {\mathcal{U(H)}},
P=\widetilde{U}^*Q\widetilde{U}\}.$$

\end{document}